\documentclass{amsart}
\usepackage{amsmath,amssymb}
\usepackage[latin1]{inputenc}
\newtheorem{theorem}{Theorem}[section]
\newtheorem{lemma}[theorem]{Lemma}

\newcommand{\N}{{\mathbb N}}
\newcommand{\Z}{{\mathbb Z}}
\newcommand{\Q}{{\mathbb Q}}

\newcommand{\U}{{\mathcal U}}

\begin{document}

\title[Hypercentral
unit groups and Hyperbolic unit groups]{Hypercentral unit groups and
the Hyperbolicity of a modular group algebra}
 \author{ E. Iwaki}
\address{ Instituto de Matem\'atica e Estat\'\i stica,
 Universidade de S\~ao Paulo,
Caixa Postal 66281, S\~ao Paulo, CEP 05315-970 - Brazil}
\email{iwaki@ime.usp.br}
\author{S.O. Juriaans}
\address{ Instituto de Matem\'atica e Estat\'\i stica,
 Universidade de S\~ao Paulo,
Caixa Postal 66281, S\~ao Paulo, CEP 05315-970 - Brazil}
\email{ostanley@ime.usp.br}
\thanks{Research partially
supported by 
FAPESP-Brazil \typeout{ and
CNPq-Brazil (Proc. 300652/95-0)}.} 
\thanks{
\hspace*{\parindent}\em{Keywords and Phrases.} {\rm group ring,
unit, normalizer, hypercenter}}

\begin{abstract}
We classify groups $G$ such that the unit group $\mathcal{U}_1(\Z
G)$ is hypercentral. In the second part, we classify groups ${G}$
whose modular group algebra has hyperbolic unit groups $V(KG)$.
\end{abstract}

\maketitle

\section{{\bf{Introduction}}}

We denote by $\Gamma = \mathcal{U}_1(\Z G)$ the group of units of
augmentation one of the integral group ring $\Z G$ of $G$.
$\mathcal{Z}_n(\Gamma)$ will denote the $n$-th centre of $\Gamma$ e
we define $\mathcal{Z}_{\infty}(\Gamma) = \bigcup_{n \in
\N}\mathcal{Z}_n(\Gamma)$. An element in
$\mathcal{Z}_{\infty}(\Gamma)$ is called an hypercentral unit.

In the finite group case, Arora, Hales and Passi in \cite{arhapa93}
showed that the central height of $\Gamma$ is at most $2$, that is,
$\mathcal{Z}_{\infty}(\Gamma) = \mathcal{Z}_2(\Gamma)$. Arora and
Passi in \cite{arpa93} then proved that
$\mathcal{Z}_{\infty}(\Gamma)=\mathcal{Z}(\Gamma)T$, where $T$
denotes the torsion subgroup of $\mathcal{Z}_{\infty}(\Gamma)$.
These results were extended to torsion groups by Li \cite{li98} and
Li,Parmenter \cite{lipar1}. In {\cite{lipar},\cite{lipar2}} they
presented some contributions to the problem for non-periodic groups.
In {\cite[Chapter VI]{her-habil}} Hertweck extended these results to
group rings $RG$ of periodic groups $G$ over $G-$adapted rings $R$.
From its exposition is clear that the containement of
$\mathcal{Z}_{\infty}(\Gamma)$ in the normalizer
$\mathcal{N}_{\Gamma}(G)$ is an important property. We present our
contribution to the study of the hypercentral units in \cite{iwaki}.
Among many other results in \cite{iwaki} it is proved that the
containement of $\mathcal{Z}_{\infty}(\Gamma)$ in
$\mathcal{N}_{\Gamma}(G)$ holds for an arbitrary group $G$. The
support of an hypercentral unit is investigated and it is proved
that the normal closure of the group generated by an hypercentral
unit is a polycyclic-by-finite group (in case, $G$ is finitely
generated).

In \cite{polcino}, Polcino Milies classified finite groups such that
the unit group of an integral group ring is nilpotent. This result
was extended to arbitrary groups by Sehgal-Zassenhaus in
\cite{sehgal}. Since nilpotent groups are hypercentral it is natural
to consider the question of classify groups $G$ such that the group
of units of an integral group ring $\mathcal{U}_1(\Z G)$ is
hypercentral. This problem was posed by several leading experts in
the field. In the second section we completely solve it as a natural
consequence of our research about the hypercentral units of an
integral group ring done in \cite{iwaki}.

In sections $3,4,5$ we deal with the topic of hyperbolic unit
groups. In the context of hyperbolic unit groups, Juriaans, Passi,
Prasad in \cite{passi} studied the groups $\mathcal{G}$ whose unit
group $\mathcal{U}(\Z \mathcal{G})$ is hyperbolic, classified the
torsion subgroups of $\mathcal{U}(\Z \mathcal{G})$ and the
polycyclic-by-finite subgroups of $\mathcal{G}$. We consider the
natural question of classifying groups $G$ for which the group of
units with augmentation one of a modular group algebra, $V(KG)$, is
hyperbolic.

Notation is mostly standard and the reader is referred to
\cite{pas,seh2} for general results on group rings. For the theory
of hyperbolic groups, we refer the reader to the reference
\cite{gromov}.

\section{Groups with Hypercentral Unit Group}
Unless otherwise stated explicitly $G$ will always denote an
arbitrary group $G$.

Firstly we recall a result proved in \cite{iwaki} which we will need
in our investigations. It will also appear in \cite{plag}.

\begin{lemma}\label{lemainv}
Let $u \in \mathcal{Z}_n(\Gamma)$ and $v$ an element of finite order
in $\Gamma$. If $c = [u,v] \not = 1$ then ${u}^{-1}vu = v^{-1}, v^2
\in G \cap \mathcal{Z}_{n+1}(\Gamma) \subseteq \mathcal{Z}_{n+1}(G),
o(v)=2^m, m \leq n, v^{2^{n-1}}$ is central and if $n=2$ then $m=2$.
In particular, elements of $\Gamma$ that are of finite order and
whose order is not a power of $2$ commute with
$\mathcal{Z}_{\infty}(\Gamma)$, and
${\mathcal{Z}}^2_{\infty}(\Gamma) \subseteq
\mathcal{C}_{\Gamma}(T(G))$, where $\mathcal{C}_{\Gamma}(T(G))$
denotes the centralizer of $T(G)$ in $\Gamma$ and $T(G)$ denotes the
set of torsion elements of $G$.
\end{lemma}

We need the following result proved in the context of nilpotent unit
groups by Sehgal-Zassenhaus in \cite{sehgal}.

\begin{lemma}\label{lemasehgal}
Suppose that $\Gamma$ is hypercentral and let $t, t_1, t_2 \in T
=T(G), g \in G$.
\begin{enumerate}
 \item Every finite subgroup of $G$ is normal in G. \label{4}
 \item If $g^{-1}tg \not = t$ then $g^{-1}tg = t^{-1}$; \label{1}
 \item If $t$ has odd order $ $then $gt=tg$;\label{2}
 \item If $1 \not = t_1$ has odd order, $t_2$ has even order
 then $T$ is a central subgroup of $G$.\label{3}
\end{enumerate}
\end{lemma}

\begin{proof}\mbox{}

$(\ref{4})$ Observe initially that since $\Gamma$ is hypercentral,
we have that $G$ is hypercentral. Since $G \subseteq \Gamma =
\mathcal{Z}_{\infty}(\Gamma)$ it follows, by Lemma
 \ref{lemainv}, that $g^{-1}tg \in \langle t \rangle$, for all $g \in
 G, t \in T$. Since every subgroup of $T$ is normal, $T$ is an abelian subgroup
 or an Hamiltonian subgroup.

 $(\ref{1})$ and $(\ref{2})$ follows immediately from Lemma \ref{lemainv}.

 $(\ref{3})$ Suppose $g \in G$ such that $g^{-1}t_1g = t_1$ and
 $g^{-1}t_2g = t^{-1}_2$. It follows that $g^{-1}t_1t_2g=
 t_1t^{-1}_2$ which should be equal to $t_1t_2$ or $(t_1t_2)^{-1}$.
 This implies that $t_1t^{-1}_2 = t_1t_2$ and $t_2^2=1$. Hence
 $g^{-1}t_2g = t_2$. Finishing the proof.
\end{proof}

We now state the main result of this section.

\begin{theorem}\label{teohyper}
 $\Gamma = \mathcal{U}_1(\Z G)$ is hypercentral if and only if $G$
 is hypercentral and the torsion subgroup $T$ of $G$ satisfies one
 of the following conditions:

 \renewcommand{\labelenumi}{(\alph{enumi})}
 \begin{enumerate}
  \item $T$ is central in $G$,\label{11}
  \item $T$ is an abelian $2-$group and for $g \in G, t \in T$
  $$
    g^{-1}tg = t^{\delta(g)}, \delta(g) = \pm{1}.
  $$\label{21}
  \item $T = K_8 \times E_2$, where $K_8$ denotes the quaternion
  group of order $8$, $E_2$ is an elementary abelian $2-$group.
  Moreover, $E_2$ is central and conjugation by $g \in G$ induces on
  $K_8$ one of the four inner automorphisms.\label{31}
 \end{enumerate}
\end{theorem}
\begin{proof}\mbox{}

 $(\Rightarrow)$ Suppose that $\Gamma = \mathcal{Z}_{\infty}(\Gamma)$. By Lemma
 \ref{lemasehgal} every subgroup of $T$ is normal in $G$ and $T$ is an abelian subgroup
 or a Hamiltonian subgroup.

 Suppose that $\Gamma$ is hypercentral and $T$ is not central. In
 any case $T$ is abelian or a Hamiltonian group.

 Suppose firstly that $T$ is a non-central Hamiltonian group with an
 element $x$ of odd order. Consider the subgroup $H = K_8 \times \langle x
 \rangle$. In this case it follows by item $(4)$ of Lemma \ref{lemasehgal} that
 $H$ is a central subgroup of $G$. Contradiction.

Suppose that $T$ is an abelian, non-central subgroup of $G$ and that
$g^{-1}t_1g = t_1$, $g^{-1}t_2g = t^{-1}_2$ for some $t_1, t_2 \in
T$. Then $g^{-1}t_1t_2g = t_1t^{-1}_2$. But by Lemma \ref{lemainv},
$g^{-1}t_1t_2g$ must be equal to either $t_1t_2$ or
$t^{-1}_1t^{-1}_2$. Hence either $t^2_1=1$ and $g^{-1}t_ig =
t^{-1}_i$, $i=1,2$ or $t^2_2 =1$ and $g^{-1}t_ig = t_i$. In any case
we obtain that $g^{-1}t_1g = t^{\delta(g)}_1$ and $g^{-1}t_2g =
t^{\delta(g)}_2$, $\delta(g) = \pm{1}$. Also $T(G)$ is an abelian
$2-$group by Lemma \ref{lemasehgal}.

Denote by $K_8 = \langle i, j \mid i^2 = j^2 = u, u^2 = 1, ji=iju
\rangle$. Every $g \in G$ maps every subgroup of $K_8$ onto itself
and induces the identity in $K_8/\langle i^2 \rangle$. In fact only
one of the four inner automorphisms
\begin{enumerate}
 \item $i \to i$, $j \to j$, $ij \to ij$.
 \item $i \to i$, $j \to ju$, $ij \to iju$.
 \item $i \to iu$, $j \to ju$, $ij \to ij$.
 \item $i \to iu$, $j \to j$, $ij \to iju$.
\end{enumerate}arises.

\mbox{}

 $(\Leftarrow)$ Under the hypothesis of the Theorem we must
prove that $\Gamma$ is hypercentral. Since $G/T$ is an ordered group
and $\Q T$ have no nilpotent elements, it follows by {\cite[Theorem
45.7]{seh2}} that $\Gamma = \mathcal{U}_1(\Z T)G$.

 We must consider three cases separately.

 \begin{enumerate}
  \item Suppose (a) holds. In this case, since $\mathcal{U}_1(\Z
  T)$ is central and $G$ is hypercentral, the result follows.
  \item Suppose (b) holds. We claim that $\mathcal{U}_1(\Z T)
  \subseteq \mathcal{Z}_{\infty}(\Gamma)$.

  $$
  \mbox{Let} \quad u \in \mathcal{U}_1(\Z{T}),\quad v=\tau{x} \in \U, \quad \tau \in
    \mathcal{U}_1(\Z{T}),\quad x \in G
  $$

Then

$$
   [u,v] =[u,\tau{x}]=[u,x]=u^{-1}u^{x} =: \gamma, \quad \mbox{where}
   \quad \gamma =1\quad \mbox{or} \quad u^{-1}u^{*}.
  $$
To see this let $u = \sum_{g \in G}\alpha_gg \in
\mathcal{U}_1(\Z{T}), x \in G.$ So $u^x = \sum_{g \in
G}\alpha_gx^{-1}gx$. Since $x$ centralizes the elements of $T$ or
$x$ acts by inversion on the elements of $T$ we obtain in the first
case that $u^x=u, \gamma = 1$ and in the second case we obtain that
$u^x=u^{*}, \gamma = u^{-1}u^{*}$. Since $T$ is abelian, it follows
that ${\gamma}^{*} = {\gamma}^{-1}$ and by {\cite[Proposition
1.3]{seh2}} $\gamma = \pm t$, for some $t \in T$. Since $\gamma$ has
augmentation $1$, $\gamma \in T$. We conclude that
$$
[\mathcal{U}_1(\Z{T}), \Gamma] \subseteq T.
$$
So
$$
[\mathcal{U}_1(\Z{T}),\Gamma, \Gamma] \subseteq [T, \Gamma] = [T,G].
$$and
$$
 [\mathcal{U}_1(\Z{T}),\Gamma, \Gamma, \Gamma] \subseteq [T,G,G].
$$Continuing this process and using the fact that $G$ is a
hypercentral group, we conclude that $\mathcal{U}_1(\Z{T}) \subseteq
\mathcal{Z}_{\infty}(\Gamma)$ and $\Gamma$ is hypercentral.
\item Suppose $(c)$. In this case, since $T$ is an Hamiltonian
$2-$group. It is well known that in this case $\mathcal{U}_1(\Z{T})$
has only trivial units. It follows that $\mathcal{U}_1(\Z{T}) = T$.
Consequently, $\Gamma = G$ is hypercentral.
 \end{enumerate}
\end{proof}


\section{Modular Group Algebras with Hyperbolic $V(KG)$}

Let ${\Z}^2$ denote the free Abelian group of rank two, $p$
be a rational prime, $GF(p^n)$ will denote the Galois Field with
$p^n$ elements, $tr.deg(K)$ denotes the transcendence degree of the
field $K$ over $GF(p)$, $\mathcal{U}(KG)$ denotes the group of units
of $KG$, $\mathcal{U}(KG)$ denotes the group of units of $KG$, $V(KG)$ denotes the group of units of $KG$ with augmentation one.

\begin{lemma}\label{lemaz2}
 Let $G$ be an arbitrary group, $K$ a field with $char(K) = p > 0$ and
 $tr.deg(K) \geq 1$. Suppose that $g_0$ is a torsion element of $G$
 and $p \nmid o(g)$. Then ${\Z}^2$ embeds in $V(KG)$ and consequently, $V(KG)$ is not hyperbolic.
\end{lemma}

In what follow we investigate under which conditions the group of
units of a modular group algebra of an arbitrary (non-trivial) group
$G$ is hyperbolic.

We denote by $\mathcal{J}(KG)$ the Jacobson Radical of $KG$ and
$\omega(G)$ represents the augmentation ideal of $KG$.

\begin{lemma}\label{lematr1}
 Suppose that $G$ is a finite (non-trivial) group and $K$ is a field,
 $char(K) = p >0$, $tr.deg(K) \geq 1$. Then
 $V(KG)$ is not hyperbolic.
\end{lemma}
\begin{theorem}\label{teofinito}
 Let $G$ be a finite (non-trivial) group, $K$ be a field, $char(K)=
 p>0$. Under these conditions, $V(KG)$ is hyperbolic if and
 only if $K$ is a finite field.
\end{theorem}

\begin{theorem}\label{teo3} Let $G$ be an arbitrary group with
torsion, $K$ be a field,
 $char(K) = p>0$. If $V(KG)$ is hyperbolic then $K$ is
 algebraic over $GF(p)$.
\end{theorem}

Our next Theorem considers the case in which $G$ is an arbitrary
(non-trivial) group, $K$ a field of $char(K)=p >0$ under the
hypothesis that $\mathcal{U}(KG)$ is
 hyperbolic.

\begin{theorem}\label{teogeral}
 Let $G$ be an arbitrary (non-trivial) group, $K$ a field of $char(K)=p >0$.
  If $\mathcal{U}(KG)$ is hyperbolic then $K$ is
 finite.
\end{theorem}

Acknowledgements: This work is part of the first authors Ph.D
thesis. He would like to thank his thesis supervisor, Prof. Dr.
Stanley Orlando Juriaans, for his guidance during this work.


\begin{thebibliography}{99}
\itemsep=-2pt
\bibitem{arhapa93}S. R. Arora, A. W. Hales,   I. B. S.  Passi, {\it Jordan decomposition and hypercentral units in integral group rings},
Comm. Algebra \textbf{21}(1993), no.1,25-35.
\bibitem{arpa93}S. R. Arora, I. B. S. Passi, { \it Central height of the unit group of an integral group ring},
Comm. Algebra \textbf{21}(1993), no.10, 3673-3683.
\bibitem{bov1} A. A. Bovdi, {\it The periodic normal divisors of the multiplicative group of a group ring I
}, Sibirsk Mat. Z. \textbf{9},(1968), no.3, 495-498.
\bibitem{bov2} A. A. Bovdi, {\it The periodic normal divisors of the multiplicative group of a group ring
II}, Sibirsk Mat. Z. \textbf{11},(1970), no.3, 492-511.
\bibitem{gromov} Gromov, M.: {\it Hyperbolic groups}, In: Essays in
group theory (S. M. Gersten, Ed.), Springer Verlag, MSRI Publ.
\textbf{8}, 1997, 75-263. MR 89e:20070.
\bibitem{her-habil} M. Hertweck, {\it Contributions to the integral
representation theory of groups},
\textsf{(http://elib.uni-stuttgart.de/opus)},2003,
Habilitationsschrift (autographed copy donated by the author).
\bibitem{plag} Martin Hertweck, E. Iwaki, E. Jespers and S. O. Juriaans, {\it On hypercentral units of integral group rings}, 2006, submitted.
\bibitem{iwaki} E. Iwaki, {\it Unidades hipercentrais em anéis de
grupo inteiro e a hiperbolicidade do grupo de unidades de uma
álgebra de grupo modular}, Ph.D. Thesis, IME-USP, 2006.
\bibitem{passi}S. O. Juriaans, I. B. S. Passi, D. Prasad., {\it
Hyperbolic unit groups}, Proc. A.M.S \textbf{133}, no. 2,
(2005),415-423.
\bibitem{li98}Y. Li, {\it The hypercentre and the n-centre of the
unit group of an integral group ring}, Canad. J. Math
\textbf{50}(1998),no. 2, 401-411.
\bibitem{lipar1}Y. Li, M. M. Parmenter, {\it Hypercentral units in
integral group rings}, Proc. Amer. Math. Soc. \textbf{129}(2001),
no. 8, 2235-2238 (electronic).
\bibitem{lipar}Y. Li, M. M. Parmenter, {\it Some results on
hypercentral units in integral group rings}, Comm. Algebra,
\textbf{31}(2003),no.7, 3207-3217.
\bibitem{lipar2}Y. Li, M. M. Parmenter, {\it The upper central series
of the unit group of an integral group ring}, Comm. Algebra,
\textbf{33}(2005), 1409-1415.
\bibitem{polcino} C. P. Milies, {\it Integral group rings with
nilpotent unit groups}, Canad. J. Math., 28:954-960, 1976.
\bibitem{pas}D. S. Passman, {\it The algebraic structure of group rings},
Robert E. Krieger Publishing Company, Malabar, Florida. Orig. Ed
1977, Reprint Ed. 1985 with corrections and appendix.
\bibitem{seh2}S. K. Sehgal, {\it Units in integral group rings}. Longman
Scientific \& Tecnical, Harlow, 1993, With an appendix by A. Weiss.
\bibitem{sehgal}S. K. Sehgal, H. Zassenhaus., {\it Integral group
rings with nilpotent unit groups}, Comm. Alg.5: 101-111, 1977.
\end{thebibliography}
\end{document}